\documentclass{elsart}

\usepackage{natbib}
\usepackage{yjsco}
\usepackage{amsfonts}
\usepackage{amscd}
\newcommand{\baseRing}[1]{\ensuremath{\mathbb{#1}}}
\newcommand{\Z}{\baseRing{Z}}
\newcommand{\R}{\baseRing{R}}
\newcommand{\C}{\baseRing{C}}
\newcommand{\N}{\baseRing{N}}
\newcommand{\Q}{\baseRing{Q}}

\newcommand{\V}{\baseRing{V}}
\newcommand{\CP}{\baseRing{P}}
\newcommand{\cstar}{\C^*}
\newcommand{\Script}[1]{\ensuremath{{\mathcal{#1}}}}

\newcommand{\ScS}{\Script{S}}

\newcommand{\NN}{\Script{N}}

\newcommand{\FF}{\Script{F}}
\newcommand{\GG}{\Script{G}}
\newcommand{\II}{\Script{I}}
\newcommand{\LL}{\Script{L}}

\newcommand{\MM}{\Script{M}}

\newcommand{\s}{\sigma}
\newcommand{\rk}{\ensuremath{{\rm rank}}}

\begin{document}
\begin{frontmatter}

\title{Restriction of $A$-Discriminants and Dual Defect Toric Varieties}
\author{Raymond Curran \thanksref{label1}}
\ead{rcurran@mscd.edu}
\thanks[label1]{Some of the results in this paper are contained
in the first author's PhD Dissertation submitted to the University of Massachusetts,
Amherst.}

\address{Department of Mathematical and Computer Sciences.
Metropolitan State College of Denver.
Denver, CO  80202, USA.}

\author{Eduardo Cattani\thanksref{label2}}
\ead{cattani@math.umass.edu}
\thanks[label2]{Partially supported by NSF Grant
DMS--0099707.  Some of the work on this paper was
done while visiting the University of Buenos Aires
supported by a Fulbright Fellowship for Research and
Lecturing. The hospitality of the Departamento de Matem\'atica, Facultad de
Ciencias Exactas y Naturales, is gratefully acknowledged.}
\address{Department of Mathematics
and Statistics. University
of Massachusetts. Amherst, MA 01003, USA}

\begin{abstract}
We study the $A$-discriminant of toric varieties.
We reduce its computation to the case of irreducible
configurations and describe its behavior under specialization
of some of the variables to zero.  We give characterizations
of dual defect toric varieties in terms of their Gale dual and classify 
dual defect toric varieties  of codimension less than or equal to four.
\end{abstract}


\end{frontmatter}

\section{Introduction}
In this paper we will study properties of the {\em sparse} or 
$A$-discriminant.  Given a configuration
$A=\{a_1,\dots,a_n\}$ of $n$ points in $\Z^d$ we may
construct an ideal $I_A \subset \C[x_1,\dots,x_n]$ and,
if $I_A$ is homogeneous, a projective toric variety $X_A \subset \CP^{n-1}$.
The dual variety $X_A^*$ is, by definition, the 
Zariski closure of the locus of
hyperplanes in $(\CP^{n-1})^*$ which are tangent to $X_A$ at
a smooth point.  Generically, $X_A^*$ is a hypersurface and
its defining equation $D_A(x)$, suitably normalized, is called 
the $A$-discriminant.  If $X_A^*$ has codimension greater than one
then $X_A$ is called a dual defect variety and we define $D_A =1$. 

The $A$-discriminant generalizes the
classical notion of the discriminant of univariate polynomials.
It was introduced by Gel'fand, Kapranov, and Zelevinsky
(their book \citep{gkz} serves as the basic reference of
our work) and it
arises naturally in a variety of contexts including 
the study of hypergeometric functions \citep{gkz2, rhf, jac} and in 
some recent formulations of mirror duality \citep{bm}.

When studying the $A$-discriminant it is often convenient
to consider a Gale dual of $A$.  This is a configuration 
$B = \{b_1,\dots,b_n\}\subset\Z^m$, where $m$ is the codimension of $X_A$ in 
$\CP^{n-1}$.  The configuration $B$, and by extension 
$A$, is said to be irreducible if
no two vectors in $B$ lie on the same line.
Equivalently, if the matroid
$\MM_B = (B,\II)$
  defined by the family,
  $\II$,  of linearly
  independent subsets of $B$ is {\em simple}.   
    In Theorem~\ref{resformula}, we prove
a univariate resultant formula which 
reduces the computation of the $A$-discriminant to the case
of irreducible configurations.  This implies, 
in particular, that
the Newton polytope of the discriminant is unchanged,
up to affine isomorphism, 
if we replace $B$ by the configuration obtained by 
adding up all subsets of collinear vectors.  This generalizes
a result of   \cite{codimtwo} for
codimension-two configurations.  We point out that, in their
case, this is a consequence of a complete description of
the Newton polytope of the discriminant.

In the study of rational hypergeometric functions,
one is interested in understanding the behavior of the
$A$-discriminant when specializing a variable $x_j$ to zero and
its relation to the discriminant of the configuration obtained
by removing the corresponding point $a_j$ from $A$.  
Theorem~\ref{specialization} generalizes the known results in
this direction (\citet[Lemma~3.2]{rhf}; \citet[Lemma~3.2]{jac}).
This specialization result was first proved by the first author in his PhD dissertation
\citep{thesis}, 
using the theory of 
coherent polyhedral subdivisions.  
We give a greatly simplified proof in \S 4, where we derive the specialization
theorem as a corollary of our resultant formula. 

Using tropical geometry methods, Dickenstein, Feitchner, and Sturmfels have
been able to compute the dimension of the dual of a projective toric variety $X_A$ and this, in particular, makes it possible to decide if a given toric variety is dual defect, i.e. if the dual variety has codimension greater than one.  Their formula \citep[Corollary~4.5]{dfs} involves the configuration $A$ and the geometric lattice, $\ScS(A)$, 
 whose elements are the 
  supports, ordered by inclusion, of the vectors in ${\rm ker}(A)$.  The
  information contained in $\ScS(A)$ is essentially the same as that
  contained in a family of flats in $\MM_B$, for a Gale dual configuration
  $B$ of $A$.
  Thus, one could say that the formula by Dickenstein, Feitchner, and Sturmfels involves both $A$ and $B$ information.   In Theorem~\ref{goodflags}, we
  use Theorem~\ref{specialization} to show that we can decide whether a configuration is dual defect  purely in terms of certain
  {\em non-splitting} flags of flats in
  the matroid $\MM_B$.  In Theorem~\ref{th:decomposition} we obtain
  a decomposition of the Gale dual configuration of a toric variety
  and give, in terms of this decomposition, a sufficient condition for the variety to be dual defect.  Although we believe this condition to also
  be necessary, we are not able to prove it at this point.
  
Dual defect varieties have been extensively studied: \citet{bel, dirocco, ein2, ein1, lanteri}. 
In particular, Dickenstein and Sturmfels have classified codimension-two dual defect varieties \citep{codimtwo} and, by completely different methods,
\citet{dirocco} has classified dual defect projective
embeddings of smooth toric varieties in terms of their
associated polytopes.  We give a complete classification of
of dual defect toric varieties of codimension less than
or equal to four in terms of the Gale duals. This implies, in particular,
that in these cases the condition in Theorem~\ref{th:decomposition} is
necessary and sufficient.
We conclude \S5 by comparing Di Rocco's list, for codimension less
than or equal to four, with our classification.

\section{Preliminaries}

We begin by setting up the notation to
be used throughout.  We will denote by $A$ a $d\times n$ integer
matrix or, equivalently, the configuration 
$A=\{a_1,\dots,a_n\}$ of $n$ points in $\Z^d$ defined by the columns
of $A$. We will always assume that $A$ has rank $d$ and set
$m := n - d$, the {\em codimension } of $A$.
Viewing $A$ as a map $\Z^n \to \Z^d$ we denote by 
$\LL_A \subset \Z^n$ the kernel of $A$. 
$\LL_A$ is a lattice of rank $m$.
For any $u\in \Z^n$ we write $u = u_+ - u_-$, where 
$u_+,u_-\in \N^n$ have disjoint support.  Let 
$I_A \subset \C[x_1,\dots,x_n]$ be the
lattice ideal defined by $\LL_A$, that is the
ideal in $\C[x_1,\dots,x_n]$ generated by all
binomials of the form:
$x^{u_+}-x^{u_-}$, where $u \in \LL_A$.  Note that for any vector
$w \in \Q^d$ in the $\Q$-rowspan of $A$ we have
$$ \langle w, u_+\rangle \ =\ \langle w, u_-\rangle$$
for all $u\in \LL_A$ and, hence, $I_A$ is $w$-weighted homogeneous.

\begin{defn} We will say that $A$ is {\em homogeneous} or
nonconfluent if the vector $(1,\dots,1)$ is in the $\Q$-rowspan of 
$A$.  
\end{defn}

Note that in terms of the configuration in $\Z^d$, $A$ is homogeneous
if and only if all the points lie in a rational hyperplane not containing the
origin.  Throughout this paper we will be interested in properties
of homogeneous configurations $A$ which depend only on the $\Q$-rowspan of $A$.  Thus, in those cases we may assume without loss of generality that
the first row of $A$ is $(1,\dots,1)$.  We shall then say that
$A$ is in {\em standard form}.

Given a homogeneous configuration $A$, let
 $X_A := \V(I_A) \subset \CP^{n-1}$ be the projective
 (though not necessarily normal) variety defined by the
 homogeneous ideal $I_A$.  The map
 $$t\in (\C^*)^d \mapsto (t^{a_1} \colon \cdots \colon t^{a_d})\in 
 X_A \subset \CP^{n-1}$$
 defines a torus embedding which makes $X_A$ into a toric variety of
 dimension $d-1$.    Generically, its
 {\em dual variety} $X_A^*$
 is an irreducible hypersurface defined over $\Z$.  Its normalized 
 defining polynomial $D_A(x_1,\dots,x_n)$ is called the sparse or
 $A$-discriminant.  It is well-defined up to sign.  If the dual
 variety $X_A^*$ has codimension greater than one, then we define
 $D_A =1$ and refer to $X_A$ as a {\em dual defect variety}
 and to $A$ as a {\em dual defect configuration}.
 Note that $X_A$, and consequently $X_A^*$, depend only on
 the rowspan of $A$.  Indeed, it is shown in
\cite[Proposition~1.2, Chapter 5]{gkz} that $X_A$ depends
only on the affine geometry of the set $A \subset \Z^d$.

Alternatively, given a configuration $A = \{a_1,\dots,a_n\}$ we consider
the generic Laurent polynomial supported on $A$:
\begin{equation}\label{eq:generic}
f_A(x;t)\ :=\ \sum_{i=1}^n x_i t^{a_i}\,, 
\end{equation}
which, for a choice of coefficients $x_i\in \C$, we view
as a regular funcion on the torus $(\C^*)^d$.  
Then, the discriminant is an irreducible polynomial in 
$\C[x_1,\dots,x_n]$ which vanishes whenever the specialization
of $f_A$ has a multiple root in the torus; i.e. $f_A$ and all
its derivatives $\partial f_A/\partial t_i$ vanishing simultaneously
at some point in  $t\in (\C^*)^d$.  Note that when $A$ is in standard form:
\begin{equation}\label{standard}
t_1 \frac{\partial f_A}{\partial t_1} \ =\  f_A
\end{equation}
and, consequently, $f_A$ and $\partial f_A/\partial t_1$ have the same
zeroes on $(\C^*)^d$.  Let $R :=\C[x][t^{\pm 1}]$ be the ring of Laurent
polynomials in $t$ whose coefficients are polynomials  in $x$, and denote
by $J(f_A)$
the ideal in $R$  generated by $f_A$ and
its partial derivatives with respect to the $t$ variables.
Set $\V_A := \V(J(f_A)) \subset \C^n_x \times (\C^*)^d_t$.  Let
   $\nabla_A$ be  the
Zariski closure of the projection of $ \V(J(f_A))$ in $\C^n_x$, then if
$\nabla_A$ is a hypersurface,
$\nabla_A = \{x: D_A(x)=0\}$.  If $A$ is homogeneous and 
$X_A$ is not dual defect then $\nabla_A$ is the cone over
$X_A^*$.

We recall that if $\nu_1,\dots,\nu_m \in \Z^n$ are a $\Z$-basis
of $\LL_A$, then the $n\times m$ matrix $B$, whose columns are
$\nu_1,\dots,\nu_m$ is called a {\em Gale dual} of $A$.  The same
name is used to denote the configuration $\{b_1,\dots,b_n\}\subset \Z^m$
of row vectors of $B$.  Gale duals are defined up to $GL(m,\Z)$-action.
We will also consider $n\times m$ integer matrices $C$, whose columns
$\xi_1,\dots,\xi_m \in \Z^n$
are a $\Q$-basis of $\LL_A\otimes_{\Z}\Q$.  In that case we will say
that $C$ is a $\Q$-dual of $A$.  For any $n\times m$ integer matrix $C$
of rank $m$ we will denote by $q$ the greatest common divisor
of all maximal minors of $C$ and call it the {\em index} of $C$.
Indeed, $q$ is the index of the lattice generated by the
row vectors of $C$, $c_1,\dots,c_n$, in $\Z^m$.
An $n\times d$ integer matrix $A$ of rank $d$ is said to be
a dual configuration of $C$ if $A\cdot C = 0$.  Note that  $C$ 
is a Gale dual of $A$ if and only if it has index $1$ and that,
if $A$
is dual to 
$C$,
then $A$ is homogeneous if and only if the row vectors of $C$ 
add up to zero.
Such a  configuration $C$  will
also be called homogeneous.
 If $c_j = 0$
for some $j$, then any dual configuration $A$ is a {\em pyramid}, i.e.
all the vectors $a_i, i\not=j$ are contained in a hyperplane.  It is
easy to check that in that case $X_A$ is dual defect.

Given an $n\times m$ integer matrix $C$
of rank $m$ we will denote by $\LL_C$ the sublattice of $\Z^n$
generated by the columns of $C$ and by
$J_C\subset\C[x_1,\dots,x_n]$ the lattice ideal
defined by $\LL_C$.  If $C$ is a Gale
dual of $A$, then $\LL_C = \LL_A$ and $I_A = J_C$ is
a prime ideal.  In any case, if 
$\xi_1,\dots,\xi_m$ are the columns of $C$ and we denote
by $J_\xi$ the ideal
$$J_\xi = \langle x^{\xi_1^+} - x^{\xi_1^-},\dots, x^{\xi_m^+} - x^{\xi_m^-}\rangle, $$
then the lattice ideal $J_C$ is the saturation
$\,J_C \ =\ J_\xi \colon (x_1\cdots x_m)^\infty.$

If $C$  is homogeneous of index $q$ then 
the variety $X_C := \V(I_C) \subset \CP^{n-1}$ has $q$ irreducible components and
they are all torus translates of $X_A=\V(I_A)$, where $A$ is a dual
of $C$.  Similarly,
the dual variety $X_C^*$ is a union of
finitely many torus translates of $X_A^*$.  In particular if one of them
is a hypersurface so is the other.  In that case, we denote by $D_C \in \C[x_1,\dots,x_n]$ the defining equation suitably normalized.    Moreover, there exist $\theta^1,\dots\theta^q \in
(\cstar)^n$ such that 
\begin{equation}\label{discc}
D_C(x) \ =\ \prod_{j=1}^q D_A(\theta^j* x)\,,
\end{equation}
where $*$ denotes component-wise multiplication.
We will say that $C$ is
dual defect if and only if $A$ is dual defect.

The computation of the $A$-discriminant is well-known in the case
of codimension-one homogeneous configurations. Let $B = (b_1,\dots,b_n)^T$,
$b_i\in \Z$,
be a Gale dual of $A$.
Reordering the columns of $A$, if necessary, we may
assume without loss of generality that $b_i >0$ for $i=1,\dots,r$ and
$b_j <0$ for $r+1 \leq j \leq n$. Set 
$$ p \ = \ b_1 + \cdots + b_r \ =  \  -(b_{r+1} + \cdots + b_n)\,.$$
Then, up to an integer factor 
\begin{equation}\label{codim-one}
D_A \ =\ \prod_{j=r+1}^n |b_j|^{|b_j|} \ \prod_{i=1}^r x_i^{b_i}
\ -\  (-1)^p \ \prod_{i=1}^r b_i^{b_i} \ \prod_{j=r+1}^n x_j^{|b_j|}
\end{equation}

We recall the notion  of {\em Horn uniformization} from 
 \cite[Chapter~9]{gkz}.  Although in \citet{gkz} this is done
 only in the case of saturated lattice ideals, the generalization
 to arbitrary lattice ideals is straightforward.  Let $C=(c_{ij})$
 be an integer matrix whose rows add up to zero, the {\em Horn map}
 $h_C:  \CP^{m-1} \rightarrow (\cstar)^m$ is defined
 by the formula 
$\ 
 h_C(\zeta_1: \cdots :\zeta_m) = (\Psi_1(\zeta), \dots, \Psi_m(\zeta)) ,
$
 where
 \begin{equation} \label{Horn2}\Psi_k(\zeta_1: \cdots : \zeta_m) = \prod_{i=1}^n (c_{i1} \zeta_1 + \cdots + c_{im}\zeta_m)^{c_{ik}}.
  \end{equation}
 We also define  $T_C: (\cstar)^n \rightarrow (\cstar)^m$ by 
 $T_C(x) := (x^{\xi_1},\dots,x^{\xi_m})$, where 
 $\xi_1,\dots,\xi_m$ are the column vectors of $C$, and set
 % \begin{equation} \label{reduceddisc}
$\   \widetilde{\nabla}_{C}:= h_C(\CP^{m-1}) \subset (\cstar)^m.$
%  \end{equation}
  
  The following result is proved in  \cite[Chapter~9, Theorem~3.3a]{gkz}
  for the case of Gale duals.  Its extension to $\Q$-duals is straightforward.
 
\smallskip

\begin{thm}\label{th:gkz3.3}
Let $A\subset \Z^n$ be a homogeneous configuration and $C \in \Z^{n\times m}$
a $\Q$-dual of $A$.  Then if $X_A^*$ is a hypersurface, so is 
$ \widetilde{\nabla}_{C}$.  Moreover, 
\begin{equation}\label{pullback}
T_C^{-1}(\widetilde{\nabla}_{C}) \ =\ \nabla_C \cap (\cstar)^n.
\end{equation}\end{thm}

\section{Discriminants and Splitting Lines}

In this section we will study the effect on the $A$-discriminant
of removing from the Gale dual configuration $B$ a set of
collinear vectors which add up to zero.  We will show that this
operation preserves the dual defect property and the Newton polytope of the
discriminant. Moreover, there is a resultant
formula relating the two discriminants.  We shall assume throughout
this section that our configurations are homogeneous.

\smallskip

\begin{thm}\label{splittinglines}
Let $A$ be a configuration in $\Z^n$ which
is not a pyramid, and $B \subset \Z^m$ a
Gale dual.  Suppose we can decompose $B$ as
$$ B\  =\  C_1 \cup C_2\,,$$
where $C_1$ and $C_2$ are homogeneous configurations, $C_1$ is
of rank $m$, and $C_2$ is of rank $1$.  Let $A_1$ be a dual of
$C_1$.  Then ${\rm codim}(\nabla_A) = {\rm codim}(\nabla_{A_1}).$
In particular,
$A$ is dual defect if and only if $A_1$ is dual
defect.
\end{thm}

\begin{pf}
Let $A_2$ be a dual of $C_2$.
We may assume without loss of generality that $A_1$ and
$A_2$ are in
standard form.  We may also assume that 
$C_1 = \{b_1,\dots,b_r\}$ and $C_2 = \{b_{r+1},\dots,b_n\}$.
Since the vectors in $C_1$ span $\Z^m$ over $\Q$, 
there is a $\Z$-relation
\begin{equation}\label{relation}
\sum_{i=1}^r \gamma_i b_i \ +\ \sum_{j=r+1}^n \mu_j b_j = 0 
\quad\hbox{
with}\quad \sum_{j=r+1}^n \mu_j b_j \not=0.
\end{equation}
It is then easy to check that the matrix
\begin{equation} \label{specialform}
A\ =\   \left(
\begin{array}{c|c}
A_1 & 0 \\ \hline
0  & A_2    \\ \hline
\begin{array}{ccc}
\gamma_1 & \cdots & \gamma_r
\end{array}  & 
\begin{array}{ccc}
\mu_{r + 1} & \cdots & \mu_{n}
\end{array}
\end{array} 
\right) 
\end{equation}
is dual to $B$ and, consequently, we may assume 
 that $A$ agrees with the matrix
(\ref{specialform}).  We can write $d = d_1 + d_2 + 1$,
where: $d_1 = r - m$ and $d_2 = n - r -1$ and view
$A_1$, $A_2$ as  configurations in $\Z^{d_1}$, 
$\Z^{d_2}$, respectively.  We let $t = (t_1,\dots,t_{d_1})$,
$s = (s_1,\dots,s_{d_2})$, $x = (x_1,\dots,x_r)$, and
$y=(y_{r+1},\dots,y_n)$.  Given $u\in \C^*$, we let
$u^\gamma * x = (u^{\gamma_1} x_1,\dots,u^{\gamma_r} x_r)$.
We define $u^\mu * y$ in an analogous way.

If $A$ is as in (\ref{specialform}), 
$
f_A(x,y;t,s,u)\ =\ f_{A_1}(u^\gamma * x; t) + f_{A_2}(u^\mu * y; s) $ and, therefore, 
$$J(f_A) \ =\ \langle J(f_{A_1}(u^\gamma* x;t)), J(f_{A_2}(u^\mu* y;s)), \partial f_A/\partial u\rangle\,.$$
In particular, we get a map
$\Phi:  \V_A \rightarrow \V_{A_1}$ given by
$\ \Phi(x,y,t,s,u)\ =\ (u^\gamma * x, t)$.
We also define $\Psi: \V_A \rightarrow \cstar \times \nabla_{A}$ by
$\ \Psi(x,y,t,s,u)\ =\ (u,x,y)\,.$
Let $Z = Im(\Psi)
\subset \cstar \times \nabla_A$, and let $\Pi: \V_{A_1}  \rightarrow \nabla_{A_1}$ denote 
the natural projection.  Finally, define 
$\phi: Z \rightarrow \nabla_{A_1}$ by
$\phi(u, x,y) = u^\gamma * x$.
Then the diagram 
\begin{equation} \label{comm}
\begin{CD}
\V_A   @>{\operatorname{\Phi}}>>   \V_{A_1}\\
@V{\operatorname{\mathit{\Psi}}}VV      @VV{\operatorname{\mathit{\Pi}}}V \\
Z   @>{\operatorname{\mathit{\phi}}}>>     \nabla_{A_1}
\end{CD}
\end{equation}
commutes.  We note that $\dim Z = \dim \nabla_A$.  Indeed, 
the natural projection $p\colon Z \to \nabla_A$ has finite
fibers since, for any $(u,x,y) \in Z$, $u^\mu* y\in \nabla_{A_2}$.
But $A_2$ is a codimension-one configuration and therefore its
discriminant is given by (\ref{codim-one}).  Hence, $u$ must satisfy
an equation of the form 
$u^q \, =\,c\,  y^\alpha$,
for some $q\in \Z$, $c\in \Q$, and $\alpha \in \Z^{n-r}$.  

We now claim that the conclusion of Theorem~\ref{splittinglines}
will follow from Lemma~\ref{gensurj}, proved below,
which asserts that $\phi$ is generically surjective with
fibers of dimension $n - r$.  Indeed,  we have
$\ \dim \nabla_A = \dim Z = \dim \nabla_{A_1} + n - r\,$
and, consequently,
$${\rm codim}(\nabla_A) = n - \dim \nabla_A = r - \dim \nabla_{A_1} = 
{\rm codim}(\nabla_{A_1}) \,.$$
\nopagebreak[4] \end{pf}

Before proving the statements on generic surjectivity and fiber dimension, 
we prove an auxiliary Lemma.

\smallskip

\begin{lem} \label{divisor} 
Let $A$ be a  $d \times n$ integer matrix of rank $d$ with 
Gale dual $B = \{b_1, \dots, b_n \}$.  Let $x \in \C^n, t \in (\C^*)^d$.
Suppose that for some 
$\Theta  \in \Z^n$,
$\,\V_A \subset \left\{ f_A (\Theta \ast x; t)= 0 \right\}.$
Then $\Theta_1 b_1 + \cdots + \Theta_n b_n =0.$  
\end{lem}  
\begin{pf}  Let $t_0=(1,\dots,1) \in (\C^*)^d$.  Then, the
set  $\{x\in \C^n : (x, t_0) \in \V_A\}$ agrees with the 
conormal space of $X_A$ at the point $[1:\cdots:1]\in X_A \subset \CP^{n-1}$.
For each such $x=(x_1,\dots,x_n)$ we have, by assumption
$$\Theta_1 x_1 + \cdots + \Theta_n x_n =0.$$
Hence $\Theta$ lies in the tangent space to $X_A$ at the point
$[1:\cdots:1]$.  Since this tangent space equals the row span of $A$,
the result follows.
\end{pf}

\begin{lem}\label{gensurj}
Under the hypotheses (\ref{specialform}), the map $\phi \colon Z \to \nabla_{A_1}$ is generically surjective
with fibers of dimension $n-r$.
\end{lem}

\begin{pf}
To prove the first statement we show that
$\Phi \colon \V_A \to \V_{A_1}$ is generically surjective.
Let $(\bar x,t) \in \V_{A_1}$ and choose $(u,y)$ such that
\begin{equation}\label{fiber}
D_{A_2} (u^\mu* y ) \ =\ 0.
\end{equation}  
As noted above, for
any choice of $y\in\C^{n-r}$ there are finitely many 
possible choices of $u$ satisfying (\ref{fiber}).  We next
choose $s \in (\C^*)^{d_2}$ such that $(u^\mu* y,s)\in \V_{A_2}$.
Note that the assumption that $A_2$ is in standard form implies
that if $(u^\mu * y,s)\in \V_{A_2}$ then so does 
$(u^\mu* y,s_\lambda)$, where 
%\begin{equation}\label{indep}
$\,s_\lambda \ =\ (\lambda s_1, s_2,\dots,s_{d_2})
\,;\ \lambda\in \C^*.$
%\end{equation}
For the given choice of $u$, let  $x$ be defined by
$ \,u^\gamma* x \ =\ \bar x\,.$
Therefore, $(x,t) \in \V(J(f_{A_1}(u^\gamma *x;t))$.   Thus,
it suffices to show that we can choose $\lambda\in \C^*$ such
that $(x,y,t,s_\lambda,u)$ satifies
\begin{equation}\label{lasteq}
\frac {\partial f_A}{\partial u} (x,y;t,s_\lambda,u) \ =\ 0
\end{equation}
But clearly
$$u\ \frac {\partial f_A}{\partial u} (x,y;t,s_\lambda,u) \ =\ 
f_{A_1}(\gamma*u^\gamma * x; t) +\lambda f_{A_2}(\mu* u^\mu*y; s),$$
where $\gamma* u^\gamma * x = (\gamma_1 u^{\gamma_1} x_1, \dots, 
\gamma_r u^{\gamma_r} x_r)$, and similarly for $\mu * u^\mu *y$.
Lemma~\ref{divisor} and (\ref{relation}) imply that we may 
assume without loss of generality that $(y,s,u)$ have been chosen so
that
$f_{A_2}(\mu* u^\mu* y; s)\not = 0$.  Thus, if $(\bar x, t)$
are so that 
$\,f_{A_1}(\gamma* \bar x; t) \not = 0,$
then we can certainly choose $\lambda\in \C^*$ so that (\ref{lasteq})
holds and, consequently, $\Phi$ is surjective outside the 
zero locus of $f_{A_1}(\gamma* \bar x; t)$.  Appealing once again
to Lemma~\ref{divisor} and (\ref{relation}), it follows that
this zero 
locus does not contain $\V_{A_1}$ which 
completes the proof of the first assertion.

Finally, we note that the remark after (\ref{fiber}) implies
the statement about the fiber dimension of $\phi$.
\end{pf}

Suppose now that we are under the same assumptions as in Theorem~\ref{splittinglines}.  That is,
$A$ is a configuration in $\Z^n$ which
is not a pyramid. $B \subset \Z^m$ is a
Gale dual of $A$ which may be decomposed as
$ \,B =  C_1 \cup C_2\,$,
where $C_1$ and $C_2$ are homogeneous configurations. $C_1$ is
of rank $m$, and $C_2$ is of rank $1$.  Moreover,
let $A_1$ be a dual of $C_1$.  We then have

\medskip

\begin{thm}\label{splittinglines:np}
If $C_1$ has index $q$,  then the Newton polytope $\NN(D_{A})$ is affinely 
isomorphic to $q\cdot\NN(D_{A_1})$.\end{thm}

\begin{pf}
By Theorem~\ref{splittinglines}, $D_A = 1$ if and only if $D_{A_1} = 1$,
thus we may assume $D_A  \neq 1$.
Let $B = \{b_1,\dots,b_n\} \subset \Z^m$ and suppose that that $C_1 = \{b_1,\dots,b_r\}$.  We will then show that
the  projection $\pi_r \colon \R^n \to \R^r$ on the first $r$ coordinates
maps $\NN(D_{A})$ to $q\cdot\NN(D_{A_1})$.  Since both of these polytopes
have the same dimension the result follows.

Note that since the vectors $\{b_{r+1},\dots,b_n\}$ are all collinear and
$b_{r+1}+\cdots+b_n=0$, we have, for all $k=1,\dots,m$, that the product
$$\prod_{i = r+1}^n (b_{i1} \zeta_1 + \cdots + b_{i,n-d}\zeta_{n-d})^{b_{ik}}$$
is a constant $\lambda_k\in\Q$.  Hence, the defining equations 
$F_B(z),  F_{C_1}(z)$ of $\widetilde{\nabla}_{B}, \widetilde{\nabla}_{{C_1}},$
are related through
\begin{equation}\label{eq:np}
F_B(z_1, \dots, z_{m}) = F_{C_1}(\lambda_1 z_1, \dots, \lambda_{m} z_{m}).
\end{equation}
By (\ref{pullback}), substituting 
$z_j$ by $x^{\nu_j}$, $j=1,\dots,m$, where $\nu_j$ is the the $j$-th column
vector of $B$, into $F_B(z)$ gives the discriminant $D_A(x)$ up to a
Laurent monomial
factor.  On the other hand, this same substitution in the right hand side
of (\ref{eq:np}) yields a polynomial in 
$\C[x_1,\dots,x_r]$ whose support equals that of 
$D_{C_1}$.  Hence 
$$\pi_r (\NN(D_A)) = \NN(D_{C_1})\,.$$
Since, on the other hand, (\ref{discc}) implies that 
$\NN(D_{C_1}) = q\cdot \NN(D_{A_1})$, the result follows.
\end{pf}

\begin{defn}\label{btilde} A configuration 
$B = \{b_1, \dots, b_n\}\subset \Z^m$ is called {\em irreducible}
if any two vectors in $B$ are linearly independent. 
If $A$ is dual to an irreducible configuration $B$, we shall
also call $A$ irreducible. Given a configuration
$B$ we will denote by $\tilde B$
the irreducible configuration obtained
by removing all vectors lying on splitting lines and replacing
non-splitting subsets of collinear vectors in $B$ by their sum.
\end{defn}

\medskip

\begin{rem}\label{matroid}
$\MM_B = (B,\II)$ be the matroid 
  defined by the family,
  $\II$,  of linearly
  independent subsets of $B$.   Then $B$ is irreducible if and
  only if $\MM_B$ is {\em simple}.
\end{rem}

\medskip

\begin{defn}\label{def:degenerate} Let $A\subset \Z^d$ be a configuration and $B\subset \Z^m$ a 
Gale dual.  $B$ is said to be {\em degenerate} if and only if 
${\rm rank}(\tilde B) < {\rm rank}(B)$.
\end{defn}

\medskip

The following corollary may be viewed as a generalization of the results in
~\cite[\S4]{codimtwo}.  

\medskip

\begin{cor} \label{sumline}
Let $A$ be a $d \times n,$  integer matrix of rank $d$ defining a 
homogeneous configuration.  Let $B = \{b_1, \dots, b_n\}$ be a  Gale
dual of $A$.  Let $\tilde B$ be as above.  Then $\NN(D_B)$ and 
$\NN(D_{\tilde B})$ are affinely isomorphic.
\end{cor}

\begin{pf}
Let $L_1,\dots,L_s$ denote the set of lines in $\R^m$ containing vectors
in $B$.  For each $j=1,\dots,s$,  let 
$$\sigma_j :=  \sum_{b_k \in B\cap L_j} b_k.$$
Consider the configuration 
$$C:= B \cup \{\sigma_1,-\sigma_1\}\cup \cdots \cup \{ \sigma_s,-\sigma_s\}\,.$$
Repeated applications of Theorem~\ref{splittinglines:np} gives that
$\NN(D_C) \cong \NN(D_B)$.  On the other hand we may also view $C$ as
$$ C =  \tilde B \cup C_1 \cup \cdots \cup C_s\,,$$
where $C_j = \{-\sigma_j\} \cup (B\cap L_j)$.  Theorem~\ref{splittinglines:np} then implies that 
$\NN(D_C) \cong \NN(D_{\tilde B})$.
\end{pf}

We next show that, with the notation and assumptions of 
Theorem~\ref{splittinglines}, there is a univariate resultant formula relating
the discriminants $D_A$ and $D_{A_1}$.

\smallskip

\begin{thm}\label{resformula}
Let $A$, $B$, $A_1$, $C_1$, and $C_2$ be as in Theorem~\ref{splittinglines} and
let $A_2$ be a dual of $C_2$. Assume moreover that $C_1$ consists of the
first $r$ vectors in $B$. Then, there exist integers 
$\delta_1,\delta_2, \gamma_1,\dots,\gamma_r,\mu_{r+1},\dots,\mu_n, M$ such that
 \begin{equation}\label{eq:resformula}
  M \, D_A (x) \ =\ Res_u (u^{\delta_1} \, D_{A_1}(u^{\gamma}*x'),
  u^{\delta_2} \, D_{A_2}(u^{\mu} * x'')),
  \end{equation}
where $x = (x_1,\dots,x_n)$,  $x' = (x_1,\dots,x_r)$, 
 $x'' = (x_{r+1},\dots,x_n)$, and
 $*$ denotes componentwise multiplication with $u^\gamma = (u^{\gamma_1},\dots,
u^{\gamma_r})$ and $u^\mu = (u^{\mu_{r+1}},\dots, 
u^{\mu_n})$. 
 \end{thm}

\begin{pf}
  If $D_A(x) = 1,$ then $D_{A_1}(x') = 1$ by Theorem~\ref{splittinglines} and (\ref{eq:resformula})
  is clearly true.

  Suppose $D_{A_1} \neq 1$.  Let $q$ be the index of $C_1$ and let
  $w$ be a $\Z$-generator of the one-dimensional lattice 
  $\Z\langle b_{r+1},\dots,b_n\rangle$.
  Since $B$ has index $1$, $q$ is the smallest 
  positive integer
  such that 
  $\,q\, w \in \Z\langle b_{1},\dots,b_r\rangle\,.$
  We can  find integers 
  $\gamma_1,\dots,\gamma_r,\mu_{r+1},\dots,\mu_n$ such that
 \begin{equation}\label{choiceofintegers}
 \gamma_1 b_1 + \cdots + \gamma_r b_r \ =\ q\,w \ =\ -\mu_{r+1} b_{r+1}
  -\cdots-\mu_{n} b_{n}\,.\end{equation}
  We may then assume that  $A$ is as in (\ref{specialform}) and therefore,
  since both $A_1$ and $A_2$ are in standard form, it follows from
  (\ref{standard}) that
if  $D_A(x) = 0$ then the discriminants
 $D_{A_1}(u^{\gamma}* x')$ and $D_{A_2}(u^{\mu}* x'')$ vanish simultaneously
 for some $u\in \C^*$.  Let $\delta_1,\delta_2\in \Z$ be such
 that $u^{\delta_1}\,D_{A_1}(u^{\gamma}* x')$ and $u^{\delta_2} D_{A_2}(u^{\mu}* x'')$ are
 polynomials in $u$ with non-zero constant term.  Then there exists a
 polynomial $F(x)$ such that
\begin{equation}\label{reseq}
Res_u (u^{\delta_1} \, D_{A_1}(u^{\gamma}* x'),
  u^{\delta_2} \, D_{A_2}(u^{\mu}* x''))\ =\ F(x)\,D_A(x)\,.
  \end{equation}
The proof of Theorem~\ref{splittinglines:np} implies that the degree
 of $D_A(x)$ in the variables $x'$ equals $q \, \deg(D_{A_1}(x'))$.
 On the other hand, the degree of the left-hand side of (\ref{reseq})
 is the $u$-degree of $u^{\delta_2}\,D_{A_2}(u^{\mu}* x'')$ times $\deg(D_{A_1}(x'))$.
 By definition of $w$, we can write $b_j = \beta_j \, w$,
 $\beta_j\in \Z$,  $j=r+1,\dots,n$, and therefore 
 $$ q = -\mu_{r+1}\,\beta_{r+1} - \cdots - \mu_{n}\,\beta_{n}$$
 but then it follows from the expression (\ref{codim-one}) for the
 discriminant of a codimension-one configuration that
 $$ \deg_u(u^{\delta_2} \,D_{A_2}(u^{\mu}* x'')) \ =\ q\,.$$
 Hence both sides of (\ref{eq:resformula}) have the same degree in
 the variables $x'$ and, consequently, $F(x)$ depends only on
 $x''=(x_{r+1},\dots,x_n)$.  
 
 Suppose $F(x'')$ is not constant.  We can write
 \begin{equation} \label{aux}
u^{\delta_1} \cdot D_{A_1}(u^\gamma* x') = g_l(x') u^\ell + \cdots + g_1(x') u + g_0(x').
\end{equation}
\begin{equation} \label{aux2} u^{\delta_2} \cdot D_{A_2  }(u^\mu* x''))=
u^{q} \prod_{\beta_j >0} \beta_j^{\beta_j} \prod_{\beta_j <0} x_j^{-\beta_j} - \prod_{\beta_j < 0} {\beta_j}^{-\beta_j} \prod_{\beta_j >0} x_j^{\beta_j}.
\end{equation}
Choose $a'' = (a_{r+1}, \dots, a_n)$ with $F(a'') = 0.$
Then $$Res_{u} (u^{\delta_1} \cdot D_{A_1}(u^\gamma* x'), u^{\delta_2} \cdot D_{A_2  }(u^\mu* a''))= 0$$
for all $x' = (x_1, \dots, x_r).$ This means Equations~(\ref{aux}) and (\ref{aux2}) are solvable
in $u$ for all $(x',a'').$  There are 
 at most $q$ possible values for $u$ which solve (\ref{aux2}), 
which means that (\ref{aux}) must be the zero polynomial which is
a contradiction since the monomials appearing in $g_i(x')$ are distinct monomials 
of $D_{A_1}.$  Thus $F(x)$ is a constant $M\in \Z$.
  \end{pf}
  
  \begin{rem}\label{coefficients}
  We note that there are many possible choices for
  $\delta_1,\delta_2, \gamma,\mu$ in Theorem~\ref{resformula}.  Indeed, it
  suffices that $\gamma$ and $\mu$ satisfy (\ref{choiceofintegers})
  and that $\delta_1$ and $\delta_2$ be chosen so that the products $u^{\delta_1}\,D_{A_1}(u^{\gamma}* x')$ and $u^{\delta_2} D_{A_2}(u^{\mu}* x'')$ be
 polynomials in $u$ with non-zero constant term.  In fact, if we replace
 (\ref{choiceofintegers}) by
\begin{equation}\label{anotherchoiceofintegers}
 \gamma'_1 b_1 + \cdots + \gamma'_r b_r \ =\ q'\,w \ =\ -\mu'_{r+1} b_{r+1}
  -\cdots-\mu'_{n} b_{n}\,,\end{equation}
where $q' = kq$, with $k$ a positive integer, then its effect is to make
a change of variable $u \mapsto u^k$ in the resultant and therefore
we would have:
\begin{equation}\label{modifiedresformula}
  M \, D_A (x)^k \ =\ Res_u (u^{\delta'_1} \, D_{A_1}(u^{\gamma'}* x'),
  u^{\delta'_2} \, D_{A_2}(u^{\mu'}* x''))\,
  \end{equation}
  for suitable integers $\delta'_1$, $\delta'_2$.
 \end{rem}
 
 \smallskip

The following corollary which will be needed in the next section
describes the effect on the discriminant of adding to the $B$ configuration a vector and its negative.

\smallskip

\begin{cor}\label{plusminus}
 Let $A \in \Z^{d\times n}$ be a homogeneous configuration and
let $B = \{b_1,\dots,b_n\}\subset \Z^m$, $m=n-d$, be a Gale dual.  Let $v \in \Z^m$ be a non-zero vector and let 
$$B^\sharp \ :=\ B \cup \{v,-v\}\,.$$
Let $A^\sharp$ be dual to $B^\sharp$.  Let $x=(x_1,\dots,x_n)$ and
$D_A \in \C[x]$, $D_{A^\sharp} \in \C[x;y_+,y_-]$, the discriminants
associated with $A$ and $A^\sharp$, respectively.  Then
$$D_A(x)\ =\ D_{A^{\sharp}}(x, y_+, y_-)|_{y_+ = 1, y_- = -1}.$$
\end{cor}

\begin{pf} Since $B$ has index $1$, we can write
 $$ v \ =\  \sum_{j=1}^n \gamma_j \, b_j\ ;\ \gamma_j\in \Z,$$
 and setting $\mu_{n+1} = 0$, $\mu_{n+2} = 1$, we can apply (\ref{eq:resformula}) and obtain
 $$ D_{A^{\sharp}}(x;y_+,y_-)\ =\ Res_u(u^{\delta_1} \cdot D_A(u^\gamma*x), y_+ + u\,y_-)$$
 for a suitable integer $\delta_1$.  We may specialize this resultant to 
 $y_+ = 1$, $y_-=-1$ since that does not change the $u$-degrees of the 
 polynomials involved and obtain:
 $$D_{A^{\sharp}}(x, y_+, y_-)|_{y_+ = 1, y_- = -1} 
 \ =\ Res_u(u^{\delta_1} \cdot D_A(u^\gamma*x), 1 - u) = D_A(x).$$
 \end{pf}
 
 We end this section 
with a simple example to illustrate how we can use Theorem~\ref{resformula}
and Corollary~\ref{plusminus} to reduce the computation of discriminants
to that of irreducible configurations and univariate resultants.

\medskip

\noindent {\bf Example.\,}  
We work directly on the $B$ side and consider a configuration  
$B$ consisting of seven vectors $\{b_1,\dots,b_7\}$, where  
$$b_1 = (0,1),\quad b_2 = (-3,1),\quad b_3 = (2,-3),\quad b_4 = (-1,1),$$  
$$b_5 = (1,0),\quad b_6 = (3,0),\quad b_7 = (-2,0).$$  
The last $3$ vectors lie on a line $L$ and $\sigma(L)=(2,0)$.  
As before, we set  
$$B^\sharp = B \cup \{\sigma(L),-\sigma(L)\} = C_1 \cup C_2,$$  
where $C_1 = \{b_1,b_2,b_3,b_4,\sigma(L)\}$ and 
$C_2 = \{b_5,b_6,b_7,-\sigma(L)\}$.  We let $\{x_1,\dots,x_7\}$ 
denote variables  
ssociated with $\{b_1,\dots,b_7\}$, respectively, and let  
$y_+$, $y_-$ be associated with $\sigma(L)$ and $-\sigma(L)$.  

We note that $C_1$ and $C_2$ are homogeneous configurations  
satisfying the assumptions in Theorem~\ref{resformula} and  
$\rm{index}(C_1) =1$.  Following the notation of    
Theorem~\ref{resformula} we have $w = (1,0)$ and therefore  
\begin{equation}\label{rel1}  
b_1 - b_4 \, =\, w\, =\, -(-1) b_5.  
\end{equation}  
On the other hand, using Singular \citep{singular} we compute  
\begin{eqnarray*}  
\lefteqn{D_{C_1}(x_1,x_2,x_3,x_4,y_+) \,  = \, 
256x_2^5 x_3^6 x_4 +13824 x_1 x_2^6 x_3^3 x_4^2 +
186624x_1^2 x_2^7 x_4^3-}\\
&& 432x_2^2x_3^8y_+^2
- 24224x_1x_2^3x_3^5x_4y_+^2 - 359856x_1^2x_2^4x_3^2x_4^2y_+^2 - 
432x_1x_3^7y_+^4 -\\ 
&&24696x_1^2x_2x_3^4x_4y_+^4
-1210104x_1^3x_2^2x_3x_4^2
y_+^4 -  823543x_1^4x_4^2y_+^6.     
\end{eqnarray*}  
While clearly  
$$D_{C_2}(x_5,x_6,x_7,y_-) \ =\ 8x_5x_6^3 - 27 x_7^2 y_-^2.$$  
Thus, given (\ref{rel1}), we may apply Theorem~\ref{resformula}  
with $\delta_1=0$, $\delta_2=1$ and obtain  
\begin{eqnarray*}
\lefteqn{D_{B^\sharp}(x,y_+,y_-) \, =\,  {\rm Res}_u(D_{C_1}(ux_1,x_2,x_3,u^{-1}x_4,y_+),  
uD_{C_2}(u^{-1}x_5,x_6,x_7,y_-))\,=}\\  
&& 5038848x_2^5x_3^6x_4x_7^6y_-^6-
     746496x_1x_3^7y_+^4x_5^2x_6^6x_7y_-^2      -421654016x_1^4x_4^2y_+^6x_5^3x_6^9
  -\\  
&& 2098680192x_1^2x_2^4x_3^2x_4^2y_+^2x_5x_6^3x_7^4y_-^4 -  
 42674688x_1^2x_2x_3^4x_4y_+^4x_5^2x_6^6x_7^2y_-^2   -\\
  &&  2519424x_2^2x_3^8y_+^2x_5x_6^3x_7^4y_-^4 -
  141274368x_1x_2^3x_3^5x_4y_+^2x_5^1x_6^3x_7^4y_-^4 + \\  
 &&  272097792x_1x_2^6x_3^3x_4^2x_7^6y_-^6 +
 3673320192x_1^2x_2^7x_4^3x_7^6y_-^6.   
\end{eqnarray*}  
According to Corollary~\ref{plusminus} setting $y_+=1, y_-=-1$ yields  
$D_B(x)$.  Since  $y_-$ appears only raised to
even powers, the expression
for $D_B(x)$ is obtained from that for $D_{B^\sharp}(x,y_+,y_-)$ erasing
$y_+$ and $y_-$. Finally, note that if, instead of (\ref{rel1}), we use  
the relation:  
$$\sigma(L) \, =\, 2w\, =\, -(-\sigma(L)) , $$  
then, as noted in Remark~\ref{coefficients}  
$$D_{B^\sharp}^2 (x,y_+,y_-) \,= \,{\rm Res}_u(D_{C_1}(x,uy_+),  
D_{C_2}(x,uy_-)).$$

 \section{Specialization of the $A$-discriminant} 
 
 The main result of this section is 
a specialization theorem for the $A$-discriminant generalizing Lemma~3.2 in \citep{rhf} and Lemma~3.2 in \citep{jac}.  In these references,
the lemmata in question play an important role in the study of rational hypergeometric functions.

We begin with a general result on the variable grouping in the
$A$-discriminant.

\smallskip

\begin{prop} \label{monomialgrouping}
Let $A$ be a $d \times n,$ integer matrix of rank $d$ and
$B=\{b_1,\dots,b_n\} \subset \Z^m$ a Gale dual of $A$. Let $D_A(x)$,
$x=(x_1,\dots,x_n)$, be the sparse discriminant.  
Then, if $b_k$ and $b_\ell$, $1\leq k,\ell\leq n$, are positive 
multiples of each other, $$D_{A} |_{x_k = 0}  = D_{A} |_{x_\ell = 0}.$$
\end{prop}
\begin{pf}
Define $\omega_k \in \R^n$ by $(\omega_k)_j = -\delta_{kj}$,  
$j = 1, \dots, n,$.   It is clear that the initial form
$in_{\omega_k}(D_A)$ of $D_A$ relative to the weight $\omega_k$
agrees with the restriction
$D_{A} |_{x_k = 0}$.   Thus, it suffices to show that
\begin{equation} \label{initial}
in_{\omega_k}(D_A)=in_{\omega_l}(D_A).
\end{equation}
 
 We recall \cite[Chapter 10, Theorem 1.4 a]{gkz} that the
secondary fan $\Sigma(A)$ is the normal fan  to the Newton polytope 
 $\NN(E_A)$ of the 
 principal $A-$determinant (we refer to ~\cite[Chapter 10]{gkz} for
 the definition and main properties of the
principal $A-$determinant). Then
 \begin{equation} \label{initial2}
 in_{\omega_k}(E_A) = in_{\omega_\ell}(E_A)
 \end{equation} if and only if 
 $\omega_k$ and $\omega_\ell$ are in the same relatively open 
 cone of $\Sigma(A)$.
 
 On the other hand, it follows from \cite[Lemma 4.2]{bfs}, that  the linear map
$-B^T: \R^n \to \R^{m}$
defines an isomorphism of fans between the secondary fan, $\Sigma(A)$,
 and its image, a  polytopal fan $\Script{F}$ defined on $\R^{n-d}.$
 Hence,   (\ref{initial2}) holds 
  if and only if $-B^T \cdot \omega_k$ and 
 $-B^T \cdot \omega_\ell$ are in the same relatively open cone
 of $\Script{F}.$  But $-B^T \cdot \omega_k = b_k$ is a positive 
 multiple of  $-B^T \cdot \omega_l = b_l$ by assumption, so 
 they must be in the same relatively open cone of $\Script{F}.$
 
 Since $D_A$ is a factor of $E_A$ by  \cite[Chapter 10, Theorem 1.2]{gkz}, 
the normal fan of $E_A$  refines that of $D_A$.  Then, any two 
 vectors giving the same initial form on $E_A$ give the same initial form
 on $D_A.$  This proves equation (\ref{initial}) and concludes the proof 
 of the Proposition.
 \end{pf}

As before, let $A = \{a_1,\dots,a_n\}$ be a homogeneous configuration
in $\Z^d$ which
is not a pyramid. 
For any index set $I\subset \{1, \dots, n\}$, we denote by $A(I)$ the 
subconfiguration of $A$ consisting of $\{a_i, i\in I\}$.
Let $B\subset \R^m$ be a Gale dual of $A$. 
Given a line $\Lambda \subset \R^m$, let 
$$I_\Lambda = \{j:  b_j \not \in \Lambda\} \subset \{1, \dots, n\}\,;\ J_\Lambda = \{1, \dots, n\}\setminus I_\Lambda.$$
and
$$
%\label{sum}
%{\displaystyle
\sigma(\Lambda) \ :=\ \sum_{j\in J_\Lambda} b_j\,.
$$
If $\Lambda$ is a non-splitting line, let $w$ be the $\Z$-generator
of $\Z\langle b_j; j\in J_\Lambda\rangle$ 
in the same direction as $\sigma(\Lambda)$ and, for $j\in J_\Lambda$ write
$b_j = \beta_j w$.  We set $J^+_\Lambda = \{j \in J_\Lambda, \beta_j > 0 \}$ and define $J^-_\Lambda$
accordingly. 

We may now prove the main result of this section

\smallskip

\begin{thm}\label{specialization}
Let $A$ be a homogeneous, $d \times n$ integer matrix of rank~$d$, and let $\Lambda$ be a 
non-splitting line.   Then, for 
any $j \in J^{+}_\Lambda$, 
$$D_{A(I_\Lambda)}  \hspace{6pt} {\textrm divides} \hspace{6pt} D_{A} |_{x_j = 0}.$$
\end{thm}

\begin{pf}
We may assume that $I_\Lambda = \{1,\dots,r\}$, and let us denote
by $x' = (x_1,\dots,x_r)$, $x''=(x_{r+1},\dots,x_n)$.
Let $B^\sharp = B \cup \{\sigma(\Lambda),-\sigma(\Lambda)\}$ and
$A^\sharp$ a dual of $B^\sharp$.  As we have done
before, let us denote by $y_+$, respectively $y_-$, the variable
associated with $\sigma(\Lambda)$, respectively $-\sigma(\Lambda)$.
By Corollary~\ref{plusminus}
$$D_A(x)\ =\ D_{A^{\sharp}}(x, y_+, y_-)|_{y_+ = 1, y_- = -1}.$$
On the other hand, we can write
$B^\sharp = C_1 \cup C_2$, where 
$$C_1 = \{b_1,\dots,b_r, \sigma(\Lambda)\}\ ;\quad
C_2 = \{b_{r+1},\dots,b_n, -\sigma(\Lambda)\}.$$
Let $w$ be a generator of $\Z\langle b_{r+1},\dots,b_n\rangle$
so that $\sigma(\Lambda) = cw$ with $c$ a positive integer.  Let $q$ 
be the index of $C_1$.  Then we may write:
$$q\cdot \sigma(\Lambda) \ =\  c\cdot q \cdot w \ =\  -q\cdot (-\sigma(\Lambda))\,.$$
Thus, it follows from (\ref{modifiedresformula}) that, up to constant,
$$
  (D_{A^\sharp} (x))^c \ =\ Res_u (D_{A_1}(x',u^q \cdot y_+),
   D_{A_2}(x'',u^q \cdot y_-))\,,
   $$ since we can choose $\delta'_1=\delta'_2=0$.
   Consequently
   $$(D_A(x))^c = Res_u (D_{A_1}(x',u^q \cdot y_+),
   D_{A_2}(x'',u^q \cdot y_-))|_{y_+ = 1, y_- = -1}.$$
   
   On the other hand, let $b_j = \beta_j\cdot w$, $\beta_j\in \Z$,
   $j=r+1,\dots,n$.  Then, since $-\sigma(\Lambda) = -c\cdot w$,
   $$D_{A_2}(x'',u^q\cdot y_-) = K_1\, \prod_{j\in J_\Lambda^+} x_j^{\beta_j} -
   K_2\,u^{cq}\,y_-^c \prod_{j\in J_\Lambda^-} x_j^{-\beta_j},$$
   where $K_1$ and $K_2$ are integers.  It then follows that
   we may specialize $x_j = 0$, $j\in J_\Lambda^+$, in the resultant
   since that does not change the leading term of $D_{A_2}(x'',u^q \cdot y_-)$.
   Hence, up to constants and monomials:
   $$(D_A(x))^c|_{x_j=0} \ =\ D_{A_1}(x',u^q \cdot y_+)^{cq}|_{u=0,y_+=1}
    \ =\ D_{A_1}(x',y_+)^{cq}|_{y_+=0}.$$
    But, since $\sigma(\Lambda)$ is the unique vector in the line $\Lambda$ in
    the configuration $C_1$, it follows that $A(I_\Lambda)$ is a non-facial
    circuit in $A_1$ and therefore by \cite[Lemma~3.2]{rhf},
    $D_{A(I_\Lambda)}$ divides $D_{A_1}(x',y_+)|_{y_+=0}$ and the
    result follows.
\end{pf}
 
\section{Dual Defect Varieties} 
 
 In this section we apply the specialization Theorem~\ref{specialization} 
 and recent
 results in \citep{dfs} to 
 prove, in Theorem~\ref{goodflags}, a Gale dual characterization of dual defect toric varieties.  This leads to a classification of dual defect toric
 varieties of codimension less than or equal to four.  Motivated by this
 classification, we prove that the Gale dual of a configuration may be
 decomposed as a disjoint union of non dual-defect configurations which
 are maximal in an appropriate sense.  Using this decomposition we give
 a sufficient condition for a configuration to be dual defect.  We believe
 that this condition is necessary as well.  Indeed, it follows from
 Theorems~\ref{three-classif} and \ref{four-classif}, that this is
 the case for codimension less than or equal to four.
 
 Throughout this section, we
 let $A$ be a homogeneous configuration of $n$ points in $\Z^d$ which
  is not a pyramid.  
  We assume moreover that the elements of $A$ span the lattice $\Z^d$.  
  As always, if convenient, we will view $A$ as a $d\times n$ integer matrix
  of rank $d$.
  Let $X_A$ denote the associated projective toric 
  variety and $X_A^* \subset \CP^{n-1}$ its dual variety.  
  Let $\ScS(A)$ denote the geometric lattice whose elements are the 
  supports, ordered by inclusion, of the vectors in ${\rm ker}(A)$.
 The following result is proved in \citep{dfs} using tropical geometry methods.
 
 \smallskip
 
  \begin{thm}[{\citet[Corollary~4.5]{dfs}}]\label{tropical} 
  Let $A$ be as above.  The dimension of 
  $X_A^*$ is one less than the largest rank of any matrix 
  $(A^t,\sigma_1,\dots,\sigma_{n-d-1})$, where 
  $\sigma_1,\dots,\sigma_{n-d-1}$ is a proper maximal chain in $\ScS(A)$.
  \end{thm}
  
 \smallskip

  Let $B\subset \Z^m$, $m = n-d$, be a Gale dual of $A$ and let 
  $\MM_B = (B,\II)$ be the matroid 
  defined by the family,
  $\II$,  of linearly
  independent subsets of $B$.  Given a subset $B' \subset B$, the
  {\em rank} of $B'$ is defined as the cardinality of the maximal element of
  $\II$ completely contained in $B'$.  A subset $F \subset B$ is called
  a $k$-{\em flat} if it is a maximal, rank-$k$ subset of $B$.  Clearly
  every subset $B'\subset B$ spans a subspace
  $\langle B' \rangle \subset \R^m$ whose dimension equals the rank of $B'$.  
  A subspace $W \subset \R^m$ is said to be $B$-spanned if 
  $\dim(W) = \rk (B\cap W)$.
  Given a flat $F \subset B$ we denote
 
  $$\sigma(F) \ =\ \sigma(\langle F \rangle) \ =\ \sum_{b\in F} b\,.$$
  A subset $C \subset B$ such that $\sigma(C)=0$ will be called a
  {\em homogeneous subconfiguration} (or a homogeneous flat if $C$ is
  a flat in $B$).
   
 \smallskip

 \begin{defn} 
 A $k$-flag of flats $\FF$ is a flag
  $F_0\subset F_1\subset\cdots\subset F_k$, where
  $F_j \subset B$ is a $j$-flat.  The flag is said to be {\em non-splitting}
  if and only if $\sigma(F_j) \not\in \langle{F_{j-1}}\rangle$, for all $j=1,\dots,k$.
\end{defn}
 
 \smallskip

Note that $F_0 = \emptyset$ and $\langle F_0 \rangle =\{0\}$, so we will usually drop it from the notation. 
If $\FF$ is a non-splitting flag then, for all $j=1,\dots,k$, $ \langle F_{j}\rangle$ is a $B$-spanned subspace and
$\sigma(F_j) \not= 0$.  Moreover, 
$\langle F_j\rangle$ projects to a non-splitting line in $\R^m/\langle F_{j-1}\rangle$.
Clearly, the projection of a non-splitting $k$-flag $\FF$ to
$\R^m/\langle F_{1}\rangle$ is a non-splitting $(k-1)$-flag in
the configuration defined by the projection of $B$.

The following is a characterization of dual defect
toric varieties which parallels that contained
in Theorem~\ref{tropical} although it only involves the Gale dual $B$.

\smallskip

\begin{thm}\label{goodflags} Let $A \subset \Z^d$ be as above and $B\subset \Z^m$ a
Gale dual of $A$.  Then $X_A$ is dual defect if and only if $B$ does not
have any non-splitting $(m-1)$-flags.
\end{thm}

\begin{pf} 
We prove the if direction by induction on the codimension $m$.  The result
is obviously true for $m=1$.  Assuming it to be true for configurations
of codimension $m-1$, let $B$ be a codimension $m$ configuration with
a non-splitting $(m-1)$-flag $F_1\subset \cdots\subset F_{m-1}$.
Let $\pi_1:\Z^m \to \Z^{m-1}$ denote the projection onto a rank
$m-1$ lattice complementary to $\langle F_1 \rangle \cap \Z^m$ and
let $G_j = \pi_1(F_{j+1})$.  Clearly, $G_1\subset\cdots\subset G_{m-2}$
is a non-splitting $(m-2)$-flag for $\pi_1(B_1)$, where 
$B_1 := \{b\in B: b\not \in F_1\}$.  We recall that $\pi_1(B_1)$ is a 
Gale dual for the configuration $A_1 := \{a_i\in A : b_i \in B_1\}$.
By induction hypothesis, $A_1$ is not dual defect and, by Theorem~\ref{specialization}, the discriminant $D_{A_1}$ must divide
an appropriate specialization of $D_A$.  Hence $A$ is not dual defect.

We also prove the converse by induction on the codimension $m$.
Once again, the case $m=1$ is clear.
We begin by considering the special case of  a
configuration $A$ with an irreducible Gale dual $B$.
If $A$ is not dual defect, by Theorem~\ref{tropical},
there exists a proper maximal chain in $\ScS(A)$,
$\sigma_1,\dots,\sigma_{n-d-1}$, such that the matrix
$M:=(A^t,\sigma_1,\dots,\sigma_{n-d-1})$ has rank $n-1$.
After reordering the columns of $A$, and consequently the entries of
$\sigma_j$, we may assume  that 
${\rm supp}(\sigma_j) = \{1,\dots,k_j\}$ with
$k_1 < \cdots < k_{n-d-1}.$

We claim that there exists an index $i$, $k_{n-d-2} < i \leq k_{n-d-1}$
such that the matrix $M_i$, obtained by removing the $i$-th row and the
last column of $M$, has rank $n-2$.  Indeed, if the columns of $M_i$ 
are linearly dependent then, since the corresponding columns of $M$
are independent, it follows that the basis vector $e_i$ may
be written as a linear combination of the first
$n-2$ columns of $M$. If this were true for every $i$, $k_{n-d-2} < i \leq k_{n-d-1}$, we could write the vector
$$ \sigma_{n-d-1} - \sigma_{n-d-2}\  =\  \sum_{k_{n-d-2} < j \leq k_{n-d-1}} e_j$$
as a linear combination of the first $n-2$ columns of $M$, a contradiction.

We fix now an index $i$, as above, such that ${\rm rank}(M_i) = n-2$.
Let $A'$ be configuration obtained by removing
the $i$-th column of $A$.  Notice that the vectors $\sigma'_1,\dots,\sigma'_{n-d-2}$
obtained, also, by removing the zero in the $i$-th entry from the 
corresponding $\sigma_j$, define a proper maximal chain in $\ScS(A')$.
We then have, by Theorem~\ref{tropical}, that $A'$ is not dual defect
and, therefore any Gale dual $B'$ of $A'$ must contain a non-splitting
$(m-2)$-flag $\,G_1\subset \cdots \subset G_{m-2}$.  Now, since $B$ is
irreducible, $B'$ agrees ---up to
$\Q$-linear isomorphism--- with the projection of $B$ onto
$\R^m/\langle b_i\rangle$.  Then, denoting by $V_j$ the lifting of
$\langle G_{j-1}\rangle$ to $\R^m$, $j=2,\dots,m-1$,  and setting
$$F_j \ :=\ V_j \cap B\ ;\quad j=2,\dots,m-1,$$
$F_1 = \{b_i\}$, we have that 
$F_1\subset F_2 \subset \cdots F_{m-1}$ is a non-splitting flag
of flats in $B$.

Finally, consider the general case. That is, let $A$ be a non dual-defect
configuration  whose Gale dual $B$ is not necessarily 
irreducible.  As before, let $\tilde B$ be the irreducible configuration
obtained from $B$ by replacing all subsets of collinear vectors
in $B$ by their sum.  Note that $\tilde B$ need not have index one, but
we may still consider a dual $A_1$ of $\tilde B$.  It follows from 
Corollary~\ref{sumline} that $D_{A_1} \not=1$.  Moreover, a Gale 
dual $B_1$ of $A_1$, being $\Q$-linearly isomorphic to $\tilde B$, is
irreducible.  Therefore $B_1$ has a non-splitting $(m-1)$-flag.  But then
so do $\tilde B$ and $B$.
\end{pf}

\begin{cor}\label{degenerate}
Let $A\subset \Z^d$ be a homogeneous configuration and let $B$ be a
Gale dual.  Then if $B$ is degenerate, $A$ is dual defect.
\end{cor}
 
 \begin{pf}
 If ${\rm codim}(A) =m$ but $B$ is degenerate, then ${\rm rank}(\tilde B)<m$
 and $\tilde B$ may not contain
 any non-splitting $(m-1)$-flags and, therefore, neither does $B$.
   \end{pf}

 Note that, by Theorem~\ref{goodflags}, if $A$ is not a pyramid and
 ${\rm codim}(A) =2$, then $D_A=1$ if
and only if a Gale dual $B$ has no non-splitting one-flags, i.e. if and
only if every
line is splitting or, equivalently, if $\tilde B = \emptyset$.
This classification of codimension-two dual defect toric varieties
 is contained in Corollary~4.5 of \citep{codimtwo}.  This observation
 may be generalized to the codimension-three case:

\smallskip

 \begin{thm}\label{three-classif} 
  Let $A= \{a_1,\dots,a_n\}\subset \Z^d$ be a homogeneous configuration of codimension three,
  which is not a pyramid.
Let $B\subset \Z^3$ be a Gale dual 
of $A$. Then $D_A = 1$ if and only if $B$ is degenerate.
\end{thm}

\begin{pf} By the above Corollary and Theorem~\ref{goodflags} it suffices
to show that if  $B$ is an irreducible configuration of
rank three, then $B$ has a non-splitting two-flag.
Let $b$ and $b'$ be distinct elements in $B$ and set 
$F_2$ be the two-flat containing $\{b,b'\}$.  If, $\sigma(F_2) \not=0$ then
we may assume $\sigma(F_2) \not\in \langle b \rangle$ and 
$\{b\}\subset F_2$ is a non-splitting two-flag.  
On the other hand, suppose every $B$-spanned plane  $P \subset \langle B \rangle$
satisfies $\s(P)=0$.  Then, fixing
an element $b\in B$, and denoting by $P_1,...,P_r$ the distinct 
$B$-spanned planes containing $b$ we would have that 
$0 = \sigma(B) = \sigma(P_1) + \cdots + \sigma(P_r)  - (r-1)\cdot b$.
But, we have assumed $\sigma(P_i)=0$ for all $i=1,\dots,r$.
Hence $r=1$, and this implies that $\rk(B)=2$, a contradiction.
\end{pf}

We consider now the case of codimension-four configurations:

\smallskip

\begin{thm}\label{four-classif}   
Let $A= \{a_1,\dots,a_n\}\subset \Z^d$ be a homogeneous configuration of codimension four,
  which is not a pyramid.
Let $B\subset \Z^4$ be a Gale dual 
of $A$. Then $D_A = 1$ if and only if  either $B$ is degenerate, or there exist 
planes $P, Q \subset \R^4$, such that $P \cap Q = \{0\}$, and 
every non-splitting line lies either in $P$ or in $Q$. 
\end{thm}

\begin{pf}
Let $A$ be such that $D_A =1$ and suppose $B$ is non-degenerate.
Let $\tilde B$ be the irreducible configuration as in
Definition~\ref{btilde}.  Since $B$ is non-degenerate the
vectors in $\tilde B$ span $\R^4$ and, by Corollary~\ref{sumline},
$D_A =1$ if and only if $D_{\tilde B}=1$.  Thus, we may assume
without loss of generality that $B$ is irreducible.  We note that if
$B = C_1 \cup C_2$,
where $C_1$ and $C_2$ are homogeneous configurations contained in 
complementary planes $P$ and $Q$, respectively, then $B$ may not
contain any non-splitting three-flags and, therefore, $A$ is dual defect.

In order to prove the only-if direction of 
Theorem~\ref{four-classif} we begin with two lemmas which  hold for arbitrary rank.

\smallskip

\begin{lem} \label{newniceflag} 
Let $B$ be a homogeneous
configuration of rank $m$ and let
$\Lambda\subset \langle B \rangle$ be a line.
Suppose $B$ has a non-splitting flag of rank $k$.
Then, $B$ has a non-splitting flag $\GG$ of rank $k$ such
that $\langle G_{k}\rangle \cap \Lambda = \{0\}$.
\end{lem}

\begin{pf}
We prove the result by induction on $k$, $1\leq k \leq m-1$.  
The result is obvious for
$k=1$ since $m\geq 2$ and $B$ is homogeneous which means that
the number of non-splitting one-flats in $B$ is either zero or
at least three.  
Assume it to be true for
 non-splitting flags of rank less than $k$,
and let $\FF$ be  a non-splitting flag $\FF$ in $B$ of
rank $k\geq 2$.
We can
assume that $\langle F_1 \rangle \not= \Lambda$. 
Consider the
projection $\pi(B)$ to $\langle B \rangle/\langle F_1 \rangle$.  
$\Lambda$ projects to a line 
$\bar\Lambda$ in $\langle B \rangle/\langle F_1\rangle$.
Moreover, the projection of $\FF$ defines a non-splitting 
flag of rank $k-1$ in $\pi(B)$.
By inductive hypothesis there exists a non-splitting flag
$\bar\GG$ in $\pi(B)$ of rank $k-1$ such
that $\langle \bar G_{k-1} \rangle \cap \bar\Lambda = \{0\}$.
Let $W_{j+1} \subset \langle B \rangle$ be the unique subspace of dimension 
$j+1$ containing $\langle F_1 \rangle$ and projecting onto
$\langle \bar G_{j} \rangle$, $j=1,\dots,k-1$.  Notice that by
construction $W_{k}\cap \Lambda \subset \langle F_1 \rangle$
but, since $\Lambda \cap \langle F_1 \rangle = \{0\}$ we have
$W_{k}\cap \Lambda = \{0\}$.  Setting $G_1 = F_1$, 
$G_j = W_j \cap B$ for $j=2,\dots,k$, we get the desired non-splitting
$k$-flag in $B$.
\end{pf}

\begin{lem} \label{existence} Let $A\subset \Z^d$ be a homogeneous
configuration of codimension $m$ and $B$ a Gale dual.   If $B$ is non-degenerate, then there exists a flat $F\subset B$ of rank
$m-1$  such that
$\sigma(F) \not = 0$.  Moreover, if we denote by $B_F$ the homogeneous
configuration in $\langle F \rangle$ defined by
$B_F :=F \cup \{-\sigma(H)\}$, then, if $B_F$ 
is non dual-defect,  $B$ is not dual defect.
\end{lem}

\begin{pf} 
 If every flat of rank $m-1$ is homogeneous, 
 let $s < m-1$ be the maximal rank of a non-homogeneous flat $F$ in $B$.  We have $s>0$ since
 $B$ is non-degenerate.
 Choose a flat $G$ of rank $s$ with $\sigma(G)\not= 0$ and
 let
 $\Theta_1,...,\Theta_r$ be the rank $s+1$ flats which contain $G$.  
By assumption, $\sigma(\Theta_i)=0$  for all $i=1,\dots,r$.  Then,
$$0 = \sigma(B) = \sum_{i=1}^r \sigma(\Theta_i) - (r-1) \cdot \sigma(G)= - (r-1) \cdot \sigma(G).$$
Hence $r=1$ and therefore $B$ has rank $s+1$.  Since
$s+1 < m$ this implies that $B$ is degenerate,   
 a contradiction.  

Suppose now that $B_H$ is not dual defect.  By Theorem~\ref{goodflags},
$B_H$ has a non-splitting flag $\GG$ of rank $m-2$ and, by
Lemma~\ref{newniceflag}, we may assume that $G_j \cap \langle\s(F) \rangle =\{0\}$. But then, 
$$G_1 \subset \cdots \subset G_{m-2} \subset F$$
is a non-splitting flag of rank $m-1$ in $B$.  Applying Theorem~\ref{goodflags} again we deduce that $B$ is not dual defect.
\end{pf}

\begin{cor}\label{distline2} Let $A\subset \Z^d$ be a homogeneous configuration of codimension four and suppose a Gale dual
$B \subset \R^4$ of $A$ is irreducible.  Suppose $B$ does not have any
non-splitting three-flags and let $F$ be a rank-three flat with
$\sigma(F) \not =0$.  Then $\sigma(F) \in F$ and 
the elements $\{b \in F : b \not= \sigma(F)\}$ span a 
plane  $P \subset \langle F \rangle$, with $\sigma(P)=0$.
\end{cor}

\begin{pf} Let 
$B_F$ be  as in Lemma~\ref{existence}. Since $B$ is dual defect
so is  $B_H$  and hence,
by Theorem~\ref{three-classif}, $B_F$ must be degenerate.
Since $F$ has rank three and $B$ is irreducible, this can only
happen if $\sigma(F) \in F$, so that 
$\{\sigma(F),-\sigma(F)\}$ define a splitting line.  The second assertion is then clear by Theorem~\ref{three-classif}.
\end{pf}

We now return to the proof of Theorem~\ref{four-classif}.  Because
of Corollary~\ref{sumline} and Theorem~\ref{goodflags}, it suffices
to prove that if $B\subset \R^4$ is an irreducible, non-degenerate
configuration which does not have any non-splitting three-flags, then
$B = C_1 \cup C_2$, where $C_1$ and $C_2$ are homogeneous, rank-two
configurations.

Let $F \subset B$ be a rank three flat with $\sigma(F) \not =0$.
By Corollary~\ref{distline2}, $F\cap B = C_1 \cup \sigma(F)$ and
$C_1$ is a rank-two flat with $\sigma(C_1)=0$.  Let 
$C_2 := B \backslash C_1$.  We claim that $C_2$ does not have 
any non-splitting two-flags.  Indeed, suppose $G_1 \subset G_2$ is 
a non-splitting two-flag.  Let $b \in C_1 \backslash G_2$.  Such $b$
exists since $C_1 \not = G_2$.  Then, letting $G_3$ be the smallest
three-flat containing $G_2 \cup \{b\}$, we would have that
$\,G_1 \subset G_2 \subset G_3\,$
would be a non-splitting three-flag in $B$, contradicting our assumption.
But, it is easy to see that the argument used in the proof of Theorem~\ref{three-classif} implies that since $C_2$ is irreducible
and has no non-splitting two-flags, it
must have rank two and $\sigma(C_2) =0$.  Since $B$ has rank four,
the planes $\langle C_1 \rangle$ and $\langle C_2 \rangle$ must be
complementary.
\end{pf}

Theorem~\ref{four-classif} motivates the following decomposition theorem which gives a sufficient condition for
a Gale configuration to be dual defect.

\begin{thm}\label{th:decomposition}
Let $B$ be a homogeneous, irreducible configuration of rank $m$.
Then, we can write 
\begin{equation}\label{eq:decomposition}
B = C_1 \cup \cdots \cup C_s,
\end{equation}
where the $C_i$'s are homogeneous,  disjoint, non dual-defect 
subconfigurations of $B$.  Moreover, $C_i$ is a flat 
 in 
$C_i\cup C_{i+1}\cup \cdots \cup C_s$ and the $C_i$'s are maximal
with these properties.  Moreover,
the rank of a non-splitting flag in $B$ is bounded by 
\begin{equation}\label{eq:decomposition2}
\rho = \rho(B) :=\sum_{i=1}^s {\rm rank}(C_i) - s.
\end{equation}
Hence
if $\rho \leq m -2$, $B$ is dual defect.
\end{thm}

\medskip

\begin{rem} It follows from Theorems~\ref{three-classif} and \ref{four-classif}
that the condition $\rho \geq m-1$ is a necessary and sufficient condition
for a configuration $B$, of rank at most four, to be dual defect.
We expect this to be the case in general.  This would give a complete classification of
dual defect toric varieties in terms of their Gale configuration.
\end{rem}

\begin{pf}
The following two lemmas, necessary for the proof of Theorem~\ref{th:decomposition},
may be of independent interest as well.

\medskip

\begin{lem}\label{niceflag2} 
Let $B $ be a homogeneous
non dual-defect configuration of rank $m$.  Suppose
$V\subset \langle B \rangle$ is a  $k$-dimensional subspace,
$0\leq k < m$.  
Then, $B$ has a non-splitting flag $\FF$ of rank $m-1$ such
that $\langle F_{1}\rangle \cap V = \{0\}$.
\end{lem}

\begin{pf}
We proceed by induction on $m$.  The result is clear for
$m=2$.  Assume our statement holds for configurations of
rank $m-1$.  Let $\GG$ be a non-splitting flag of rank
$m-1$ in $B$.  If $\langle G_{1}\rangle \cap V = \{0\}$ we
are done.  Assume then that $G_1 \subset V$ and consider
the projection $\pi(B)$ to  $\langle B \rangle /\langle F_1 \rangle$.
Then, $\pi(B)$ is not dual defect and, by inductive hypothesis,
there exists a non-splitting $(m-2)$-flag $\bar\FF$ in $\pi(B)$
such that $\langle\bar F_1 \rangle\cap \pi(V) = \{0\}$.  Let
$W_{j+1}$ be the unique subspace of   $\langle B \rangle$ containing
$\langle G_1\rangle$ and projecting to $\langle \bar F_j\rangle$ 
and set $F_{j+1} = W_{j+1}\cap B$.  Note that 
$\sigma(F_{j+1}) \not \in \langle  F_j\rangle$ since 
$\bar\FF$ is non-splitting.  Now $\langle\bar F_1\rangle\cap \pi(V ) = \{0\}$
implies that $\langle  F_2\rangle\cap  V  = \langle G_1 \rangle$.
Now, since $F_2$ is spanned by non-splitting one-flats, there exists
a  one-flat $F_1 \subset F_2$,  with
$\langle F_1 \rangle \not= \langle G_1\rangle$, and such that $\sigma(F_2) \not \in F_1$. The 
flag $F_1\subset F_2 \subset \cdots \subset F_{m-1}$ is a non-splitting
flag in $B$ with $\langle F_1 \rangle \cap V =\{0\}$.
\end{pf}

\begin{lem}\label{splitflat}
Let $B$ be an irreducible,  homogeneous, dual defect configuration and let $\Lambda$ a line in $\langle B \rangle$.
Then there exists a homogeneous, non dual-defect
flat $C \subset B$ of rank $k$, $2\leq k <m$, such that $\langle C \rangle \cap \Lambda =\{0\}$.
\end{lem}

\begin{pf} We proceed by induction on
$m =  {\rm rank}(B)$. If $m\leq 3$ then,
by Theorem~\ref{three-classif}, there are no irreducible, non dual-defect
configurations.  So assume that $m\geq 4$ and that  the result holds for configurations
of rank less than $m$.  Let 
$k<m-1$ be the largest rank of a non-splitting flag
in $B$.    We may assume that $k\geq 2$.  
Otherwise, given any one-flat $F_1$ in $B$, every two-flat containing it
must be homogeneous, but
this is 
impossible since $B$ is irreducible.  Moreover, by Lemma~\ref{newniceflag},
we may assume that $B$ has a non-splitting $k$-flag $\FF$ such
that $\langle F_k \rangle \cap \Lambda = \{0\}$.

Let $\Theta_0,\dots,\Theta_q$ be the distinct 
$(k+1)$-flats in $B$ containing $F_k$. Since $m > k+1$, $q\geq 1$, and at most
one $(k+1)$-flat may contain both $\langle F_k \rangle$ and $ \Lambda$.
Hence may assume
 $\Lambda \cap \langle \Theta_j\rangle=\{0\}$ for $j\geq 1$.   
 If $\sigma(\Theta_j) =0$ for some $j\geq 1$, then we can take
 $C = \Theta_j$ and we are done.  If not, 
 let Let $W = \langle \Theta_1 \rangle$ and $B_W = \Theta_1 \cup \{-\sigma(\Theta_1)\}$.  
Then $B_W$ is a homogeneous configuration of rank $k+1$,
which may or may not be irreducible.  
Let $\tilde B_W$ be as in Definition~\ref{btilde}.  

Suppose $ {\rm rank}(\tilde B_W)=k$.  Then, since $B$ is irreducible,
$C:=\tilde B_W$ is a homogeneous $B$-flat of rank $k$ which, we claim, is not dual defect.
Indeed, let $j$ be
such that
  $\sigma(\Theta_1) \in F_j\backslash F_{j-1}$, we can define a non-splitting flag
$F'_1\subset \cdots \subset F'_{k-1},$
 of rank $k-1$ in $C$, by $F'_i = F_i$ for $i<j$ and
 $F'_i = F_{i+1}\cap C$ for $i=j,\dots,k-1$.
 
 If, on the other hand, $ {\rm rank}(\tilde B_W)=k+1$, then note that
 $\tilde B_W$ is dual defect.  Indeed, suppose $\tilde B_W$ has
 a non-splitting $k$-flag
 $G_1\subset \cdots \subset G_k$.  Then, by Lemma~\ref{newniceflag}, we may assume without loss of generality that $\langle G_k\rangle \cap \langle \sigma(\Theta_1)\rangle = \{0\}$.  But then $G_1\subset \cdots \subset G_k\subset \Theta_1$ would be a non-splitting flag of rank $k+1$ in $B$,
 a contradiction. Hence, by inductive hypothesis, $\tilde B_W$ has a homogeneous,
 non dual-defect flat $C$ of rank at least two and such that
 $\langle C \rangle \cap \langle \sigma(\Theta_1)\rangle =\{0\}$.  Therefore,
 $C$ is a flat in $B$ as well and the proof is complete.
\end{pf}

\medskip

We return now to the proof of Theorem~\ref{th:decomposition}.
We prove the existence of (\ref{eq:decomposition}) by
induction on the rank $m$.  If $m=2$ then, being irreducible, 
$B$ is not dual defect and we may take $B=C_1$.  

Suppose the theorem
holds for configurations of rank less than $m$ and let $B$ be
an irreducible, dual defect configuration of rank $m$.  By  
Lemma~\ref{splitflat}, there exists a homogeneous, non dual-defect,
$B$-flat $C_1 \subset B$. We may assume that $C_1$ is not
contained in any larger, homogeneous, non dual-defect $B$-flat and
$\rk(C_1)<m$. Let $B_1 = B \backslash C_1$.  Clearly,
$B_1$ is homogeneous and irreducible.  If $B_1$ is not dual defect
then taking $C_2 = B_1$ we are done.  On the other hand,
if  $B_1$ is dual defect and of rank less than $m$, then we may
apply the inductive hypothesis to write 
$B_1 = C_2 \cup \cdots \cup C_s$ where the $C_j$ are maximal, homogeneous,
disjoint,
non dual-defect 
subconfigurations of $B_1$ and, for $i\geq 2$, 
$C_i$ is a flat in 
$C_i\cup C_{i+1}\cup \cdots C_s$.   Finally, if ${\rm rank} (B_1) = m$,
we
repeat the argument and write $B_1$ as a disjoint union $B_1 = C_2 \cup B_2$, where 
$C_2$ is a homogeneous non dual-defect $B_1$ flat.
Since at each step the cardinality of the remaining homogeneous
configuration $B_j$ strictly decreases, it is clear that this
process terminates.

In order to prove the second assertion, consider
a non-splitting flag $\FF$ of rank $k$ in $B$.  We claim that,
for each $p\leq k$,
there exist $C_i$-flats $F_{i,p}\subset C_i \cap F_p$ such that
\begin{enumerate}
\item $\langle F_p \rangle = \langle F_{1,p} \rangle \oplus \cdots \oplus
\langle F_{s,p} \rangle$ and
\item if $\sigma(F_{i,p}) \in \langle F_{i,p-1} \rangle$, then 
$F_{i,p} = F_{i,p-1}$.
\end{enumerate}
Clearly, this would imply the result since the
distinct flats among the  $F_{i,p}$, $p=1,\dots,k$ would define a
non-splitting flag in $C_i$ whose rank would, therefore, be
bounded by ${\rm rank}(C_i)-1$.  To prove
the claim we proceed by induction on $p$.  If $p=1$, then we may assume
$F_1 \subset C_1$ and it suffices to choose $F_{1,1} = F_1$ and
$F_{i,1} = \emptyset$ for $i>1$.  Suppose now that we have constructed
$F_{i,p-1}$, $i=1,\dots,s$ and set $G_{i,p} := C_i \cap F_p$.
Then $F_p$ is the disjoint union of the 
$C_i$-flats $G_{i,p}$, for $i=1,\dots,s$.
Let $i_0$ be the first index such that $\sigma(G_{i_0,p}) \not\in \langle F_{p-1}\rangle$. Such an index exists since $\FF$ is a non-splitting flag.
Since $\sigma(G_{i_0,p}) \not\in \langle F_{p-1}\rangle$, there exists
a $C_{i_0}$-flat $F_{i_0,p}$ such that 
$F_{i_0,p-1} \subset F_{i_0,p} \subset G_{i_0,p}$ and 
${\rm rank}(F_{i_0,p}) = 1+ {\rm rank}(F_{i_0,p-1})$. 
Set $F_{i,p} = F_{i,p-1}$ for $i\not= i_0$.  Note that since $F_{i_0,p} \not \subset F_{p-1}$  ${\rm rank}(F_{1,p} \cup \cdots \cup F_{s,p})$ must be strictly
larger than $p-1$.  Hence
$\langle F_p \rangle = \langle F_{1,p} \rangle + \cdots +
\langle F_{s,p} \rangle$ and, for dimensional reasons, this must be
a direct sum.
\end{pf}

We have shown in Theorem~\ref{four-classif} that  if $A\subset \Z^d$ is a dual defect homogeneous configuration of codimension four,
  which is not a pyramid, and $B$ is a Gale dual then either $B$ is degenerate
  or $\tilde B = C_1 \cup C_2$, and $\langle C_1 \rangle$ and
   $\langle C_2 \rangle$ are complementary planes.  In
   this case, if $\tilde A$ is a dual of $\tilde B$ then
   $\tilde A$ is a union of homogeneous, codimension-two configurations
   lying in complementary subspaces of $\Z^d$.  Similarly, if $B$ is a degenerate
configuration consisting of vectors in a splitting line and in a complementary
three-dimensional space, then $A$ is a union of two homogeneous configurations,
of codimension one and three respectively,  lying in complementary subspaces of $\Z^d$. In either case, the projective toric variety $X_A$ is obtained
from a join of two varieties by attaching codimension-one configurations according to  
(\ref{specialform}).

More generally, if $B$  is decomposed as in (\ref{eq:decomposition}) and
$A$ is
a dual of $B$, then $A$ will be a Cayley configuration of $s$ configurations
$A_0,\dots,A_{s-1}$ in $\Z^{q}$, where $q = |B| - {\rm rank}(B) - s$, in the
following sense:

\begin{defn} Let $A_0, \dots, A_k \subset \Z^r$ be 
configurations.  The  configuration
$$ {\rm Cay}(A_0, \dots, A_k) := (\{e_0\}\times A_0) \cup \cdots \cup 
(\{e_k\}\times A_k)\subset \Z^{k+1}\times \Z^r\,,$$
where $e_0,\dots,e_k$ is the standard basis of $\Z^{k+1}$, is
called the Cayley configuration of $A_0, \dots, A_k$.
\end{defn}

In the special case  when
$B = C_1 \cup \cdots \cup C_s$, as in Theorem~\ref{th:decomposition}, is an irreducible configuration such that
$$\langle B \rangle = \langle C_1 \rangle \oplus \cdots \oplus \langle C_s \rangle,$$
then, if $A\subset \Z^d$ is dual to $B$, the toric variety $X_A$ is a join
of varieties $X_{A_1},\dots,X_{A_s}$ lying in disjoint linear subspaces and
the dual variety $X^*_A$ has codimension $s$.  However, as the following example shows, for codimension greater than four, it is no longer true that 
every dual defect
toric variety is obtained from a join by attaching codimension-one configurations according to  
(\ref{specialform}).

\smallskip

\noindent{\bf Example.} Let $A$ be the  Cayley configuration  in $\Z^4$,
$$A \, := \, {\rm Cay}(\{0,1,2\},\{0,1,2\},\{0,1,2\}).$$ The variety
$X_A$ is a smooth
three-fold in $\CP^8$. It is easy to show that a
Gale dual $B \subset \Z^5$ may be decomposed as $B = C_1 \cup C_2 \cup C_3$, where
$C_i$ is an irreducible, homogeneous, codimension-two configuration and,
therefore, non dual-defect.  Let $\rho(B)$ be as in (\ref{eq:decomposition2}).
Then $\rho(B) = 3 = \rk(B) -2$ and, by Theorem~\ref{th:decomposition}, $B$ is dual defect.  In fact using Theorem~\ref{tropical} one can show that $X^*_A$ is
a six-dimensional subvariety of $\CP^8$.
 \smallskip

Di Rocco has obtained a 
classification of dual defect projective
embeddings of smooth toric varieties in terms of their
associated polytopes
 \citep{dirocco}.  
Recall that  a homogeneous  configuration 
$A$ is said to be
{\em saturated} if
$A = \{a_1,\dots,a_n\}$
consists of all the integer points of a $d-1$
dimensional polytope with integer vertices,
$P$, 
lying on a hyperplane off the origin. Moreover, the projective
toric variety $X_A$ is smooth, if and only if the polytope $P$ is 
{\em Delzant}, that is,
for each vertex $v$ of $P$,   
there exist $w_1, \dots, w_d \in \Z^d,$ such that $\{w_1, \dots, w_d\}$ 
is a lattice basis of $\Z^d$, and $P = v + \sum_{j = 1}^d \R_+ \cdot w_j$
near $v$.  It is well known that projective
embeddings of smooth toric varieties are in one-to-one correspondence
with Delzant polytopes.

\smallskip

Di Rocco's classification theorem  
\citep[Theorem~5.12]{dirocco},
which is proved by techniques completely different to the ones
in this paper,
may now be stated as follows:

\begin{thm} \label{dir-classif}
Let $A$ be a saturated,
homogeneous, configuration in $\Z^d$ which is not
a pyramid and such that $P = {\rm conv}(A)$ is 
Delzant.  Then $A$ is dual defect
if and only if 
$$A = {\rm Cay}(A_0, \dots, A_k),$$ where
$k$ is such that 
$\max(2, \frac{d}{2}) \leq k \leq d-1$,  $A_0, \dots, A_k$
are saturated and the polytopes  $P_i := {\rm conv}(A_i)\subset \R^{d-k-1}$ are all
Delzant polytopes of the same
combinatorial type.
\end{thm}

Thus, we see that the smoothness condition puts very strong
conditions on the type of Cayley configuration we may consider.
To illustrate this, we will
list all smooth dual defect projective
toric varieties of codimension at most four.  

We note first of
all that in these cases, the configurations $A_i$ in Theorem~\ref{dir-classif}
must be one-dimensional.  In fact, let $A$ be a dual defect, saturated,
homogeneous, configuration in $\Z^d$ which is not
a pyramid and such that $P = {\rm conv}(A)$ is 
Delzant, and write
$A = {\rm Cay}(A_0, \dots, A_k),$ as in Theorem~\ref{dir-classif}.
Then, if ${\rm codim}(X_A) \leq 5$, each polytope $P_i$ must be 
one-dimensional.  
Indeed, let us consider the simplest case when the polytopes
$P_i$ are two-dimensional.  Then $d = k + 3$ and since
by assumption $k \geq (k+3)/2$, we must have $k\geq 3$.
The fewest number of integral points in a
Delzant polytope in $\R^2$ is three.  Hence 
$n = |A| \geq 12$ and $m = n - 6 \geq 6$.

Let $[p]$ denote the configuration
$\{0,1,\dots,p\}\subset \Z$.
An easy counting argument now shows that the smooth dual defect toric varieties of codimension less than or equal to four are the ones associated with the Cayley configurations listed below:

\smallskip

\newcounter{fig}
\begin{list}{\bfseries\upshape Codimension \arabic{fig}:}
{\usecounter{fig} \setcounter{fig}{1}
\setlength{\parsep}{0.25cm} 
\setlength{\labelwidth}{2cm}
\setlength{\labelsep}{0.2cm}
\setlength{\leftmargin}{2.35cm}
}
\item ${\rm Cay}([1],[1],[1])$.
\item ${\rm Cay}([1],[1],[2])$; ${\rm Cay}([1],[1],[1],[1])$.
\item ${\rm Cay}([1],[2],[2])$;   ${\rm Cay}([1],[1],[3])$; ${\rm Cay}([1],[1],[1],[2])$;\newline ${\rm Cay}([1],[1],[1],[1],[1])$.
\end{list}

\smallskip

The Gale duals of the  configurations in the above list are easily
computed.   Indeed, it is easy to see that each Cayley factor 
$A_i=[1]$ contributes a splitting line containing two vectors
from $B$, and this vectors are primitive relative to the lattice
$\Z^m$.  Similarly, each factor $A_j=[k]$ contributes a homogeneous
subconfiguration $C_j$ of rank $k$ and containing exactly $k+1$ primitive
vectors in $B$.  Thus, for example, in the codimension four case, the
configuration ${\rm Cay}([1],[2],[2])$ has a Gale dual $B$ whose
reduced configuration $\tilde B$ decomposes as $C_1 \cup C_2$, where
$C_i$ are homogeneous configurations of rank two, lying in complementary
planes, and consisting of three primitive vectors each.

\medskip

\noindent{\bf Acknowledgements.}  We are grateful to Alicia Dickenstein and to two anonymous referees  for their thoughtful comments on a previous version
of this paper.

% \bibitem[Names(Year)]{label} or \bibitem[Names(Year)Long names]{label}.
% (\harvarditem{Name}{Year}{label} is also supported.)
% Text of bibliographic item

\def\cprime{$'$}

\end{document}